\documentclass{article}

\title{Poly-infix Operators and Operator Families}
\date{}
\usepackage{stmaryrd,amssymb,amsmath,latexsym,amsthm,url,enumitem}
\usepackage[english]{babel}
\usepackage{eurosym}

\usepackage[usenames]{color}

\newcommand{\PSI}{\mathop{\Psi}}

\author{Jan A.\ Bergstra ~\&~ Alban Ponse
\\
{\small  Section Theory of Computer Science, Informatics Institute}\\
{\small  Faculty of Science, University of Amsterdam}\\
{\small \url{https://staff.fnwi.uva.nl/{j.a.bergstra/,a.ponse/}}}
}

\begin{document}

\maketitle
\thispagestyle{empty}

\begin{abstract}
Poly-infix operators and operator families are introduced as an alternative for 
working modulo associativity and the corresponding  bracket deletion convention. 
Poly-infix operators represent the basic intuition of
repetitively connecting an ordered sequence of entities with the same connecting 
primitive.
\noindent  
\end{abstract}

\bigskip\bigskip

{ \tableofcontents}

\newpage

\section{Introduction}
The source of inspiration for writing this paper we found 
in~\cite{LugtenEtAl2015} where a 
seemingly compelling case is made for the use of sums of the form 
\[n+n+n,\quad n+n+n+n, \quad n+n+n+n+n,\quad ...
\]
and the way in which these expressions are created implies that brackets do not 
enter the picture in a meaningful way.
In~\cite{LugtenEtAl2015} only substitution instances like $7+7+7$ are used, and a
reader of this source is
supposed to know that $7+7+7$ is an expression that involves three occurrences of 7. 

We will take an interest in the syntactic analysis of expressions like $7+7+7$ from the 
perspective that, e.g.\ in this particular case, we are looking
at a single occurrence of a three-place operator
rather than at two nested occurrences of a two-place infix operator with omitted bracketing
(assuming associativity) or default bracketing (from the left or from the right).

We will call an operator that is made up from a repeated but unbracketed use of the 
same infix operator a \emph{poly-infix operator}. All such operators together for a given 
infix operator constitute a poly-infix operator \emph{family}. 
The two-place operator from which 
all operators from a poly-infix operator family are made up is called its \emph{kernel}.
The use of + as a poly-infix operator is standard in Dutch primary education, see
e.g.~\cite{Beterr,Sommenf}.
The remarkable aspect of~\cite{LugtenEtAl2015} which led us to the notion of a poly-infix operator is that it introduces expressions like $2+2+2+2$ in such a manner that an explanation 
of the structure of that expression as say $(2+2) + (2+2)$ is manifestly a detour, there simply is no role for brackets in the setting of~\cite{LugtenEtAl2015}.

\subsection{Notation and axioms  for poly-infix operator families}
Given a sort name $S$ and a function symbol $\Psi$, the poly-infix operator family with 
kernel $\Psi$ for sort $S$ contains for each positive number $n\geq 2$ an operator 
\[\Psi_n: S^n \to S.\]

One axiom and two axiom schemes \eqref{AttL} and \eqref{AttR}
(Association to the Left, and Association to the Right, respectively)
with $n\geq 2$ 
are required for poly-infix operator families:\footnote{In some cases, it can be elegant,
  practical, or natural, to also introduce $\Psi_1$ defined by $\Psi_1(x)= x$, and, if a
  unit $e_\Psi$ for $\Psi$ is available, $\Psi_0$ defined by $\Psi_0=e_\Psi$.}
\begin{align}
\nonumber
\Psi_2(x_1,x_2)&= x_1\PSI x_2,\\
\label{AttL}
\tag{AttL$_{n+1}$}
\Psi_{n+1}(x_1,x_2,...,x_{n+1})
&=
\Psi_{n}(\Psi_2(x_1,x_2),...,x_{n+1}),
\\
\label{AttR}
\tag{AttR$_{n+1}$}
\Psi_{n+1}(x_1,...,x_{n},x_{n+1})
&=
\Psi_{n}(x_1,...,\Psi_2(x_{n},x_{n+1})).
\end{align}
With induction to $n$ it can be shown that for all positive numbers $n\geq 2$,
\begin{align*}
\Psi_{n+1}(x_1,...,x_n,x_{n+1})=\Psi_{2}(\Psi_n(x_1,..,x_{n}),x_{n+1}),\\
\Psi_{n+1}(x_1,x_2,...,x_{n+1})=\Psi_{2}(x_1,\Psi_n(x_2,..,x_{n+1})).
\end{align*}

\subsection{Poly-infix notation}
As an alternative, and more appealing notation, it is plausible to use infix notation for 
$-\PSI-$ and to write 
\[\text{$x_1 \PSI x_2~ ... ~x_{n+1} \PSI x_{n+2}$ 
\quad for\quad 
$\Psi_{n+2}(x_1,x_2,..,x_{n+1},x_{n+2})$.}
\]
This is called poly-infix notation for $-\PSI-$. 

Associativity of $\Psi$ as a binary operator follows immediately from the axioms of 
poly-infix families. It is obvious that both schemes are needed for associativity, by 
considering any example with $\Psi$ a non-associative function.

Using poly-infix notation, it can be proven that the mentioned axiom schemes allow the 
introduction of arbitrary bracketing in an expression of the form 
$x_1 \PSI~ ... ~\PSI x_{n+1}$. 
For instance:
\[x_1 \PSI x_2 \PSI x_3 \PSI x_{4} \PSI x_{5}= x_1 \PSI (x_2 \PSI x_3) \PSI x_{4}\PSI x_{5}.
\]

\section{Six examples of poly-infix operator kernels}
We will call an operator \emph{pre-arithmetical} if it can be reasonably introduced 
in an incremental hierarchy of datatype descriptions as well as in the explanation of that hierarchy 
in advance of an explanation of any arithmetical datatype.
For the case of poly-infix operator kernels we found three cases where a pre-arithmetical 
introduction makes sense in our view. The property of being pre-arithmetical is an informal 
one and it depends on one's view of the relation between algebra and arithmetic. 
We assume the surrounding container view of algebra with respect to arithmetic
(see~\cite{Bergstra2015}). That view reads as follows: 
\begin{enumerate}[label=\alph*)]
\item 
Arithmetic is about structures labelled as arithmetics, 
\item
each arithmetical structure is an algebra as well,  
\item 
each arithmetical algebra is surrounded by a plurality of extended, restricted, and modified 
structures which are among the non-arithmetical algebras, 
\item 
which algebras precisely must or may be 
classified as arithmetical may depend on the individual views which may differ from 
person to person, and
\item
when taking that variation of judgement into account, and assuming that empirical work to 
determine the variation in a certain condition in a statistical manner, one finds that a 
degree of arithmetically may be assigned to an algebra rather than a sharp distinction. 
\end{enumerate}

For instance some may not count the meadow of rational complex numbers 
(see~\cite{BergstraT2007,BergstraBP2015}) as an arithmetic and others may do, which gives 
it a degree of arithmeticality between 0 and 1.


\subsection{Three pre-arithmetical poly-infix operator kernels}
Three most important examples of pre-arithmetical poly-infix operator families are these:
\begin{description}
\item [parallel composition:] 
$p_1 \mid\mid .... \mid\mid p_n$ represents the parallel composition (also called merge) of 
$n$ objects (processes, entities).\footnote{%
  For an algebraic treatment of parallel composition we refer to~\cite{BergstraP2001}.}
\item [sequential composition:]  $u_1;...;u_n$ represents the sequential composition
of $n$ instructions.\footnote{%
  A more systematic name is text sequential composition, thus indicating that sequentiality 
  in time of an effectuation process is not meant. For a theoretical account of these 
  matters see e.g.~\cite{BergstraM2012}.}
\item [frame composition:] $f_1 \oplus... \oplus f_n$ represents the combination of $n$ frames.\footnote{%
  If edges and edge labels are deleted combinations of frames are mere sets of named objects,
  i.e. nodes. Frame composition extends naive set-theoretic union with deletable features 
  for labeled directed graphs. See~\cite{BP95a,BP95b}.}
\end{description}
The importance of these examples lies in the observation that expressions involving such 
operators can be found by looking at real life scenes or objects. 
The parallel presence of a number of static entities may be expressed as a merge of atoms 
for each entity. Straight-line computer programs may be modelled with an appropriate 
application of sequential composition for instruction sequences. 
A directed graph with labeled nodes and edges can be viewed as a sum of atomic frames.

\subsection{Two intra-arithmetical poly-infix operator kernels}
An operator kernel is called \emph{intra-arithmetical} if it can be introduced on the basis of 
an initial fragment of an introduction of the arithmetic algebras, and if in addition it 
contributes to the further development of the development of a complete hierarchy of 
arithmetical algebras. Two operations stand out as examples; addition 
($-+-$) and multiplication ($-\cdot-$).

Repeated addition, say $7+7+7+7+7$ can be used to explain multiplication, which is the 
approach taken in~\cite{LugtenEtAl2015}, and repeated multiplication,
say $5 \cdot 5 \cdot 5 \cdot 5$ can be used to explain exponentiation.

\subsection{A post-arithmetical poly-infix operator kernel}
Matrix multiplication in the context of $n$ dimensional square matrices of the rational 
numbers constitutes a prominent example of a post-arithmetical poly-infix operator kernel.

\section{Poly-infix operator specific equational logic}
A non-$\Psi$ expression is either a variable, or a constant, or an expression surrounded 
with brackets, or an expression with an operator different from $\Psi$ as its leading 
function symbol. A $\Psi$-expression is an application of a poly-infix operator with kernel 
$\Psi$ to a sequence of non-$\Psi$ expressions. 
The \emph{length} of a $\Psi$-expression is the 
number of arguments to its leading occurrence of a poly-infix operator with kernel $\Psi$.
Non-$\Psi$ expressions are considered $\Psi$ expressions of length 1.

A $\Psi$-expression context $C_\Psi[-]$  is a $\Psi$-expression with a hole in it, the hole 
being denoted with $[-]$. The length of $C_\Psi[-]$ is the number of top-level occurrences 
of $\Psi$ in it plus one.
For instance,
\[x \PSI [-] \PSI y\PSI z
\]
is a $\Psi$-expression context with length 4. 
Substitution of a $\Psi$-expression $P$ of length $n$ in a context $C_\Psi[-]$ of length $m$ 
produces $C_\Psi[P]$,  a $\Psi$-expression of length $m+n-1$, for instance: 
\[x \PSI [u \PSI v] \PSI y\PSI z = x \PSI u \PSI v \PSI y\PSI z.\]
 
An important  rule of equational reasoning applies in this case:
\begin{equation}
\label{rule:1}
\textit{Let $C_\Psi[-]$ be a $\Psi$-expression context. If $P=Q$, then
$C_\Psi[P] =C_\Psi[Q]$.}
\end{equation}
The virtue of this derived rule is to allow derivations without any 
manipulation of bracket pairs. 

\subsection{A formal example}
As an application of rule~\eqref{rule:1}, which may be understood as a formal 
underpinning of the work 
in~\cite{LugtenEtAl2015} we notice that it allows to derive
\[2+2+3 = 4+3\] 
from $2+2 =4$. Furthermore, we notice that by considering repeated addition in this way, no 
talk about associativity or about conventions of not writing bracket pairs which at closer 
inspection must be assumed to be present is needed.

\subsection{Motivation of the example}
The relevance of the example above is the following: 
if a teacher is aware of the theory of poly-infix operators and operator families 
as discussed above, then:
\begin{enumerate}[label=\roman*)]
\item  (s)he can speak with confidence and correctly about repeated addition as a primitive 
operator which allows its own dedicated equational logic, 
\item this exposition strategy is meaningful too
in the direction of an audience of students who have not been exposed to that theory and who 
will not be exposed to that theory,
\item further when introducing expressions like $2+2+2+2+2$ (s)he may assume with full 
confidence that in terms of talk of syntax and expressions one is dealing with a five-place 
operator, which is a member of the poly-infix operator family with $-+-$ as its kernel,
\item and as a consequence (s)he may be in no doubt that $2+2+2+2+2$ is an expression, in 
particular a $+$-expression of length five, and finally
\item an equation like $2+2+2+2+2 = 2+2 + 4 +2$ can be understood as stating of two 
expressions that these have the same value, while a useful and dedicated proof rule for 
deriving that equation is readily available.
\end{enumerate}

\section{Concluding remarks}
The use of a dedicated equational logic given the reading of an occurrence of unbracketed 
repetition of an infix operator as an occurrence of a poly-infix operator falls within the 
degrees of freedom allowed by the relativism of~\cite{Shapiro2014}. 
Pragmatically speaking, logic does not come for free given an algebra or an arithmetic, 
it calls for dedicated design, even if in a more fundamental sense there are few degrees 
of freedom.

Reasoning about poly-infix expressions, while thinking of these as values, introduces 
complications that are comparable with the resolution of inconsistencies that arise if an 
expression-oriented view on fractions and a value-oriented view on fractions are mixed 
naively. In~\cite{BergstraB2015} the latter topic is dealt with by means of an application of 
paraconsistent reasoning in the style of the chunk and permeate paradigm of~\cite{BrownP2004}.

The motivation for this work has come from three sides and can be characterized as follows: 
\begin{enumerate}[label=\alph*)]
\item 
to provide a theoretic foundation of 
the educational proposals made in~\cite{LugtenEtAl2015} which to the best of our knowledge 
necessitates considering poly-infix operators as first class citizens, 
\item
to proceed with 
efforts on contrasting an expression-oriented view and a value-oriented view on topics in 
elementary mathematics such as reported in~\cite{BergstraB2015} by taking an 
expression-oriented view on repeated occurrences of the same infix operator seriously, 
and 
\item
to provide foundations for the efforts described in~\cite{Jeene2015} on defining new 
and ``low'' reference levels of arithmetical competence which are expected to be helpful 
for improving the organization of special education and remedial teaching in the Netherlands.
\end{enumerate} 
It seems that the proposals in~\cite{Jeene2015} are entirely consistent with the approach 
taken in~\cite{LugtenEtAl2015} regarding the objective to avoid the use of brackets in
connection with what we have proposed to call ``poly-infix operators''.

\bigskip

We end with a philosophical remark. 
Albert Visser from Utrecht University has explained to one of us (JAB)
that the use of poly-infix 
operators introduces strong typing in a way which on the long run may prove 
counter-productive in view of its complexity. In particular when modelling natural language 
strong typing may lead to impractical complexity. Flexible arity allows an operator to have 
variable number and structure of arguments. Personally, I prefer strong typing, where each 
operator has a fixed number of arguments, over flexible typing because of its clarity, but 
admittedly that is a point of view which may need to be compromised once the proximity to 
natural language increases.\footnote{Indeed, poly-infix operators are implicitly used 
  in the process of learning the written form of a natural language based on a Roman alphabet,
  in the construction of words from letters, of sentences from words 
  (and punctuation symbols with special spacing 
  directives), and the construction of texts from sentences. According to some
  conventions, Arabic numerals can also be used as words in sentences (and as letters in words 
  in sms-messaging (text messaging)).}
Albert Visser has also pointed out that Leibniz has contemplated the representation of natural 
numbers with poly-infix addition applied to units, including the remarkably simple definition 
of addition in that context. Frege amongst others has criticized that view noticing that more 
precision with brackets is needed.

\end{document}